\documentclass[12pt,reqno,intlimits]{amsart}

\usepackage{amssymb} 

\usepackage{amssymb, amsmath, amsfonts}
\usepackage{bm}
\usepackage[all,2cell]{xy}\UseAllTwocells

\let\trr=\triangleright

\def\genfd{{\bf k}} 
\def\id{{\rm id}}
\def\Id{{\rm Id}}
\def\AA{{\mathcal A}}
\def\BB{{\mathcal B}}
\def\CC{{\mathcal C}}
\def\DD{{\mathcal D}}
\def\GG{{\mathcal G}}
\def\MM{{\mathcal M}}
\def\NN{{\mathcal N}}
\def\PP{{\mathcal P}}

\def\End{\operatorname{\mathrm{End}}}

\def\tto{\Rightarrow}

\def\qcoh{{\mathfrak{Qcoh}}}

\def\cat{\mathbf{Cat}}

\def\catset{\mathbf{Set}}

\def\nxpoint{\refstepcounter{subsection}%
  \makepoint{\thesubsection}}
\def\nxsubpoint{\refstepcounter{subsubsection}%
  \makepoint{\thesubsubsection}}
\def\refpoint#1{{\rm\textbf{\ref{#1}}}}

\def\makepoint#1{\medbreak\noindent{\bf #1. }}
\newcommand{\nodo}[1]{}

\begin{document}

\title[Some Equivariant Constructions]{Some Equivariant Constructions in Noncommutative Algebraic Geometry}

\author[Z. \v{S}koda]{Zoran \v{S}koda}

\maketitle

\hfill {\it To Mamuka Jibladze on occasion of his 50th birthday}

\begin{abstract}
We here present rudiments of an approach to geometric actions in
noncommutative algebraic geometry, based on geometrically
admissible actions of monoidal categories.  This generalizes the
usual (co)module algebras over Hopf algebras which provide affine
examples. We introduce a compatibility of monoidal actions and
localizations which is a distributive law. There are satisfactory
notions of equivariant objects, noncommutative fiber bundles and
quotients in this setup.
\bigskip



\end{abstract}

\bigskip

Modern mathematics is ever creating new kinds of geometries, and
again viewpoints of unification emerge. Somehow category
theory seems to be very effective in making this new order.
Grothendieck taught us how important for geometry the {\it relative}
method is and the emphasis on general maps rather than just on the
incidence hierarchy of subspaces, intersections and so on.
Important properties of maps are often just categorical
properties of morphisms in a category (possibly with a structure).
Various spectral constructions in algebra and category theory use
valuations, ideals, special kinds of modules, coreflective
subcategories, and so on, to single out genuine ``underlying
sets'' of points, or of subschemes -- to rings, algebras and
categories -- which appear as objects representing 'spaces'.
Abstract localization enables us to consider local properties of
objects in categorical setup; sheaf theory and generalizations
enable passage between local and global. It is always enjoyable
listening about a rich vision of categorical and any other
geometry from Mamuka, due to his enthusiasm, width of interests
and knowledge.

\section{Noncommutative Algebraic Geometry}\label{sec:ncg}

\nxpoint Descriptively, a {\it noncommutative space} $X$ is a
geometric entity which is determined by a structure (algebra
$A_X$, category $C_X\dots$) carried by the collection of objects
(functions, cocycles, modules, sheaves~$\dots$) which are
heuristically, or in a genuine model, living over $X$. In this
article, our primary interest will be spaces represented by
abelian categories ``of quasicoherent sheaves''.
Gabriel--Rosenberg theorem says that every scheme can be
reconstructed up to isomorphism of schemes from its category of
quasicoherent sheaves. This involves spectral constructions
\cite{Ros:lecs}: from an abelian category, Rosenberg constructs
a genuine set, its spectrum (many different spectra have been
defined for various purposes), which can be equipped with a
natural induced topology and a stack of local categories.

\nxpoint Noncommutative analogues of group actions, quotients and
principal bundles have been abundantly studied earlier,
particularly within quantum group renaissance
\cite{Drinf:ICM,Majid}, in the context of study of noncommutative
algebras and graded algebras representing noncommutative affine or
projective varieties. As known from commutative geometry, it is
easy to get out of these categories when performing the most basic
constructions, e.g., the quotient spaces. The Tannakian
reconstruction points out the correspondence between group-like
objects and categories of representations, and it is natural to
try to extend this principle not only to symmetry objects but also
to actions themselves, considering thus the actions of monoidal
categories of modules over symmetry objects to some other
categories of quasicoherent sheaves. However, not every action
qualifies.

\nxpoint (Affine morphisms.) Given a ring $R$, denote by $R-{\rm
Mod}$ the category of left $R$-modules. To a morphism of rings
$f^\circ : R \to S$ (which is thought of as a dual morphism to
$f:\mathrm{Spec}\,S\to\mathrm{Spec}\,R$) one associates
\begin{itemize}
\item {\bf extension of scalars} \index{extension of scalars} $f^*
: R-{\rm Mod} \to S-{\rm Mod}$, $M \mapsto S \otimes_R M$;

\item {\bf restriction of scalars} \index{restriction of scalars}
(forgetful functor)
$f_* :  S-{\rm Mod} \to R-{\rm Mod}$, ${}_S M \mapsto {}_R M$;

\item $f^! : R-{\rm Mod} \to S-{\rm Mod}$,
$M \mapsto {\rm Hom}_R ({}_R S, M)$.

\end{itemize}
Denote $F \dashv G$ when a functor $F$ is left adjoint to a functor
$G$. An easy fact: $f^* \dashv f_* \dashv f^!$. In particular, $f^*$
is left exact, $f^!$ right exact and $f_*$ exact. Moreover, $f_*$
is faithful. As maps of commutative rings correspond precisely to
maps of affine schemes, one says that an (additive) functor $f^*$
is {\bf almost affine} if it has a right adjoint $f_*$ which is
faithful and that $f^*$ is {\bf affine} if, in addition, $f_*$ has
a right adjoint as well (another motivation for this definition:
Serre's affinity criterion, \'{E}l\'{e}ments de g\'{e}om\'{e}trie alg\'{e}brique, II 5.2.1, IV 1.7.18).

\nxpoint (Pseudogeometry of functors) \label{s:pseudo} Given two
abelian categories ${\mathcal A}$, ${\mathcal B}$, (equivalent to small
categories) a {\bf morphism} $f : {\mathcal B} \rightarrow {\mathcal A}$
(viewed as a categorical analogue of a map of spectra or rings) is
an isomorphism class of right exact additive functors from ${\mathcal
A}$ to ${\mathcal B}$. An {\bf inverse image functor} $f^*: {\mathcal
A}\rightarrow {\mathcal B}$ of $f$ is a chosen representative of $f$.
If it has a right adjoint, then it will be referred to the {\bf
direct image functor} of the morphism $f$. An inverse image functor
$f^*$ is said to be flat (resp. coflat; biflat) if it has a right
adjoint and it is exact (resp. if $f_*$ is exact; if both $f^*$
and $f_*$ are exact). A morphism is flat (resp. coflat, biflat,
almost affine, affine) if its inverse image functor is such.

\nxsubpoint Grothendieck topologies and their noncommutative
generalizations (\cite{Ros:NcSS,KontsRos}) may be used to talk
about {\it locally affine} noncommutative spaces. For instance,
localization functors and exactness properties of functors may be
used to define relative noncommutative schemes, see
\refpoint{s:ncscheme} and \cite{Ros:NcSch}.

\section{Actegories, Biactegories, Distributive Laws}

\nxpoint To fix the notation we recall that a {\bf monoidal
category}  is given by a 6-tuple $\tilde{\mathcal C} = ({\mathcal
C},\otimes, {\bf 1}, a, \rho, \lambda)$, where $\mathcal C$ is a
category, $\otimes : \CC\times\CC\to\CC$ the monoidal product,
with unit object $\bf 1$, associativity coherence $a:\_
\otimes(\_\otimes\_)\Rightarrow (\_\otimes\_)\otimes\_$ and $\rho
: {\rm Id}_{\mathcal C} \Rightarrow {\rm Id}_{\mathcal C}\otimes {\bf 1}$
and $\lambda : {\rm Id}_{\mathcal C} \Rightarrow {\bf 1} \otimes {\rm
Id}_{\mathcal C}$ are the right and left unit coherences satisfying
the usual coherence diagrams. A (strong) {\bf monoidal functor}
$F:(\CC,\otimes,{\bf 1},a,\lambda,\rho)\!\to \!(\DD,\otimes',{\bf
1}',a',\lambda',\rho')$ is given by a triple $(F,\chi,\xi)$ where
$F:{\mathcal C}\!\to\! {\mathcal D}$ is a functor, and
$\chi:F(\_)\otimes' F(\_)\Rightarrow F(\_\otimes\_)$ and $\xi :
F(\_)\otimes'{\bf 1}'\Rightarrow F(\_\otimes {\bf 1})$ are invertible natural
transformations satisfying the coherence conditions as
in~\cite{MacLane}. An action of monoidal category $\tilde\CC$ on a
category $\AA$ is a monoidal functor from $\tilde\CC$ to the
strict monoidal category of endofunctors $\End\AA$; these data are
also said to form a (left) $\tilde\CC$-{\bf actegory}. Right actegories
correspond to reversing the order of the tensor product in
$\End\AA$. If $\mathcal L : \CC\to\End\AA$ is an action, then one
often describes it in terms of the bifunctor
$\triangleright:\CC\times\AA\to\CC$ given by $c\,\triangleright a
\!=\! \mathcal L(c)(a)$. Then the coherences 
$\chi^{\mathcal L},\xi^{\mathcal L}$ for $\mathcal L$ 
are replaced by the coherences $\Psi,u$ with components
$\Psi_M^{X,Y} : (X\otimes Y)\triangleright M \!\mapsto\!
X\triangleright (Y\triangleright M)$ and $u_M : M\to
M\triangleright {\bf 1}$ for $\triangleright$. 
Thus a $\tilde\CC$-actegory can be described as a 4-tuple 
$(X,\triangleright,\Psi,u)$.

\nxpoint \label{s:restriction} \!(Restriction for actegories). Let
$(J,\zeta,\xi)\!:\!(\BB,\otimes,\!\bm 1, a, l,r)\!\to\! (\GG,\otimes'\!,\!\bm
1', a', l',r')$ be a monoidal functor, then the bifunctor
$\triangleright_\BB = J\otimes'\Id_\GG : \BB\times\GG\to\GG$ is an
action making $\GG$ into a left $\BB$-actegory with obvious
coherence. More generally, let $\NN$ be a left $\GG$-actegory.
Then it becomes a left $\BB$-actegory as follows. The action
functor is of course $\BB\trr_\BB\NN := J\trr_\GG\NN
:\BB\times\NN\to\NN$. The action coherence component
$\Psi^\BB_{b,b',n}$ is the composition
\begin{gather*}
b\triangleright_\BB(b'\triangleright_\BB n) = Jb\trr_{\GG}(Jb'\trr_\GG
n) \stackrel{\Psi^\GG_{Jb,Jb',n}}\longrightarrow (Jb\otimes'
Jb')\trr_\GG n \stackrel{\zeta_{b,b'}\otimes n}\longrightarrow
J(b\otimes b')\trr_\GG n \\
= (b\otimes b')\trr_\BB n
\end{gather*}
We say that $\NN$ carries the {\bf restricted action} of $\BB$
along $(J,\zeta,\xi)$; we obtain $\BB$-actegory
$J_*(\NN)=(J,\zeta,\xi)_*(\NN)$ (sometimes written simply
${}_{\BB}\NN$). It is easy to check that any $\GG$-equivariant
functor $(K,\gamma) : \NN\to\PP$ of $\GG$-actegories restricts to
the $\BB$-equivariant functor $J_*(K,\gamma) :=
(K,\gamma^J):{}_{\BB}\NN\to{}_{\BB}\PP$, where the restricted
coherence $\gamma^J : \BB\trr_\BB K\to K\circ(\BB\trr_\BB\NN)$ has
components $(\gamma^J)^b_n = \gamma^{Jb}_{n} : b\trr_\BB Kn =
Jb\trr_\GG Kn \to K(Jb\trr_\GG n) = K(b\trr_\BB n)$.

\nxpoint Given two monoidal categories, $\tilde\CC$ and
$\tilde\DD$ acting on the same category $\AA$, from the left and
right via bifunctors $\triangleleft$ and $\triangleright$,
respectively, a distributive law~is a transformation $ l \!:\!
\CC\triangleright(\AA\triangleleft\DD)\!\tto\!
(\CC\triangleright\AA)\!\triangleleft \DD $ satisfying two
co\-he\-rence pentagons and two triangles, generalizing the
coherences for the usual distributive laws between monads
(\cite{Semtriples,Borceux,skoda:ncloclec}). For clarity, we draw
one of the pentagons:
$$\xymatrix{
(c\otimes c')\triangleright (a\triangleleft
d)\ar[r]^{\Psi}\ar[d]^l&
c\triangleright(c'\triangleright(a\triangleleft
d))\ar[r]^{c\triangleright l}& c\triangleright ((c'\triangleright
a)\triangleleft d)\ar[d]^l
\\
((c\otimes c')\triangleright a)\triangleleft
d)\ar[rr]^{\Psi\triangleleft d}&&(c\triangleright
(c'\triangleright a))\triangleleft d
 }$$
commutes for all objects $a\in\AA$, $c,c'\in\CC$, $d\in \DD$
(subscripts on $l,\Psi$ omitted).

We say that such data form a {\bf
$\tilde\CC$-$\tilde\DD$-biactegory} if the components of the
distributive law involved are invertible. Biactegories are a
categorification of bimodules (over monoids).

\nxpoint There exists a hierarchy in generality: distributive laws
for two actions of two different monoidal categories are more
general than between a monad and an actegory \cite{Skoda:distr},
which is in turn more general than between two monads, all
provided we allow not only strong actions but also general (lax)
monoidal and (colax) (co)monoidal functors into $\End\AA$ (e.g., a
monad interpreted as a lax monoidal functor from the trivial
monoidal category $\tilde 1$ or alternatively, as the action of a
PRO for monoids on $\AA$). For fixed $\tilde\CC$,
$\tilde\CC$-actegories, colax $\tilde\CC$-equivariant functors and
transformations of equivariant functors form a 2-category
$\tilde\CC-\mathrm{act}_c$, and a monad in that 2-category (in the
sense of formal theory of monads \cite{Street:formal}) is
precisely the usual monad in $\AA$ equipped with the distributive
law between the action of $\tilde\CC$ and the monads. The
Eilenberg-Moore construction exists for such $\CC$-equivariant
monads; this existence is an abstract consequence of a theorem on
limits for lax morphisms in \cite{Lack:laxlimits}; we have given a
direct proof and the concrete formulas for the Eilenberg-Moore
2-isomorphism from a 2-category of $\CC$-equivariant monads to
$\tilde\CC-\mathrm{act}_c$ in \cite{Skoda:eqdis}.

\nxsubpoint (Remark.) One sometimes needs more general
2-categorical symmetry objects than monoidal categories and
bicategories; hence the distributive laws between the {\it
actions} of two such 2-symmetries (each given by a pseudomonad) on
the same object may be of interest. To this aim we sketch in
\cite{Skoda:eqdis} a new concept of {\it relative distributive
law} (which is of course different form the notion of a
distributive law between two pseudomonads).

\nxpoint (Tensor product of actegories) It is the basic
observation in our work on {\it biactegories} (in preparation) 
that one can define a
tensor product of a left $\tilde\CC$-actegory $\MM$ and a right
$\tilde\CC$-actegory $\NN$ as the vertex
$\MM\otimes_{\tilde{\CC}}\NN$ of the pseudocoequalizer
$(\MM\otimes_{\tilde{\CC}}\NN,p,\sigma)$ of the functors
$\triangleright\times\NN,\MM\times\triangleleft:
\MM\times\CC\times\NN\to\MM\times\NN$ in $\cat$ with projection
$p:\MM\times\NN\to\MM\otimes\NN$ and the invertible 2-cell part
$\sigma:p\circ(\triangleright\times\NN)\tto
p\circ(\MM\times\triangleleft)$. It suffices actually to consider
the pseudocoequalizer as a representative of a bicoequalizer. If
one or both actegories underly biactegories, then the unused
actions get inherited by this tensor product; if both, a
distributive law will be induced as well, yielding therefore the
tensor product of biactegories.

For the tensor product of biactegories it is essential that we
require that the distributive laws in the definition of
biactegories be indeed invertible.

\nxsubpoint {\bf Claim.} {\it The monoidal categories,
biactegories, colax biequivariant functors, natural
transformations of colax biequivariant functors with the tensor
product of biactegories $($using pseudocoequalizers$)$ form
a tricategory $\mathrm{Biact}_c$. } \vskip .1in

This tricategory is an analogue of the bicategory of rings and
bimodules. Though straightforward, the proof is extremely long,
and left out for a future article. The result of course generalizes
to actions of bicategories instead of monoidal categories. If the
monoidal category is in fact a categorical group $G$, then one may want
to restrict to biactegories whose left and right actegories are
in fact $G$-2-torsors; this way we get bi-2-torsors, and a tricategory
relevant to the categorified geometrical Morita theory (I.~Bakovi\'c
has studied 2-torsors with structure bigroupoid in~\cite{Bakovic} and
is now thinking further on bi-2-torsors).

\nxpoint The most important special case of the tensor product of
biactegories is the construction of the {\it induction for
actegories}, which supplies the left pseudoadjoint to the
restriction 2-functor from \refpoint{s:restriction} (in the setup
of functors between appropriate 2-categories of actegories over
fixed monoidal categories). In the theory of categorified bundles
this sort of induction may be used to define the associated
2-vector bundles to the 2-torsors over categorical groups.

In the setting of~\refpoint{s:restriction}, and with $\MM$ a right
$\tilde\BB$-actegory, the pseudocoequalizer
$$\xymatrix{
\MM\times\BB\times\GG \ar@<-.5ex>[rr]_{\MM\times\triangleright}
\ar@<.5ex>[rr]^{\triangleleft\times\GG} &&\MM\times\GG\ar[rr]^-{p}
&&\MM\otimes_\BB \GG = \mathrm{Ind}^\GG_\BB\MM}$$ with
$\sigma:p\circ(\MM\times\triangleright)\tto
p\circ(\triangleleft\times\GG)$ is equipped with a canonical right
$\tilde\GG$-action, defining the induction 2-functor. It takes a
considerable work to prove the coherence pentagon for the induced
$\tilde\GG$-action.

\nxpoint {\bf Proposition.} {\it Every biactegory in
$\mathrm{Biact}_c(\tilde{\GG},\tilde{\GG}')$ is biequivalent to a
biactegory with the identity as a distributive law. }
\vskip .07in

Indeed, one replaces $\MM$ by $\GG\otimes_\GG\MM$. After
consideration of the standard construction of pseudocoequlizers in
$\cat$, one easily realizes that for the distributive law on
$\GG\otimes_\GG\MM$ one should choose identity.  
$\MM$ and $\GG\otimes_\GG\MM$ are, of
course, biequivalent biactegories.

\section{Actions of Monoidal Categories in Noncommutative Geometry}

\nxpoint (Free versus tensor product). In commutative algebraic
geometry, the category of commutative Hopf algebras over a field
$\genfd$ is antiequivalent to the category of affine group
$\genfd$-schemes, and essentially the only non-affine examples of
group $\genfd$-schemes are abelian varieties. Thus extending the
view that the affine noncommutative schemes make a category
$\mathrm{NAff}_\genfd$ dual to the category of noncommutative
rings, V. Drinfeld in the 1980-s took the viewpoint that the
noncommutative Hopf algebras are the (duals to) affine group
schemes in the noncommutative world
(\cite{Drinf:ICM,Majid,ParshallWang}). This viewpoint seemed very
successful in view of many examples, many of which are called
quantum groups. However, the drawback of this point of view is
that for many geometric constructions the product of 'spaces'
which behaves well for the development of advanced constructions
is the categorical product. In $\mathrm{NAff}_\genfd$, the latter
corresponds to the coproduct of noncommutative $\genfd$-algebras,
in other words the free product $\star$ of $\genfd$-algebras,
hence {\it not} the tensor product $\otimes_\genfd$. Examples of
cogroup objects in the category of $\genfd$-algebras exist as
well; the prime example is the noncommutative $GL_n$,
corepresenting the functor $R\mapsto GL_n(R)$ which to
$\genfd$-algebra $R$ assigns the set of all $n\times n$ invertible
matrices over $R$. However such examples are obviously very big,
close to the free algebras
(\cite{Berstein:cogroups,Fresse:cogroups}), and far from the
``quantum'' examples which are deformations of and hence closer in
size and ring-theoretic properties to commutative algebras.

\nxpoint (Bimodules as morphisms). Another peculiarity of
noncommutative geometry, the Morita equivalence, comes partly at
rescue for Hopf algebras. Indeed, geometrically and physically
meaningful constructions usually do not distinguish algebras in
the same Morita equivalence class. Thus one can compose a usual
morphism of rings with a Morita equivalence and still have a valid
morphism in the noncommutative world. In other words, one
considers bimodules as morphisms, and more generally, allowing for
nonaffine schemes, the pairs of adjoint functors between
'categories of quasicoherent sheaves'. Working over a fixed base
category (typically: category of modules over a possibly noncommutative ring
$\genfd$) sometimes restores a distinguished element in the Morita
equivalence class, namely the inverse image functor of the
morphism to the base scheme, applied to the distinguished
generator in the base. 

In this setting of 'spaces' represented by
categories over a fixed category $\mathrm{Spec}\,\genfd$, and with
adjoint pairs as morphisms, one reintroduces the Hopf
$\genfd$-algebra $H$ (where $\genfd$ is commutative) in the
disguise of the monoidal category ${}_H\MM$ of left $H$-modules
equipped with the inverse and direct image functors of a morphism
to ${}_\genfd\MM$; the direct image functor is the forgetful
functor; it is crucial that this functor is {\it strict monoidal}.
A distinguished action of the monoidal category ${}_H\MM$ on
${}_\genfd\MM$, is given by applying the direct image functor in
the first component and then tensoring in ${}_\genfd\MM$. This is
natural because the Hopf algebra $H$ lives in ${}_\genfd\MM$, and
the actions have to respect the $\genfd$-structure. Thus if
${}_H\MM$ is acting on any other category $\mathcal C$ over
${}_\genfd\MM$  the square
$$\xymatrix{
{}_H\MM \times {\mathcal C}\ar[rr]^-{\triangleright}\ar[d]
&& {\mathcal C}\ar[d]
\\
{}_H\MM\times {}_\genfd\MM\ar[rr]^-{\triangleright_0} && {}_\genfd\MM }$$
commutes,
where $\triangleright_0$ is the distinguished action. This
picture where the monoidal action represents the action of group
schemes may also be found in commutative geometry, where one
essentially takes the action of $\qcoh_G$ on $\qcoh_X$, where $G$
is a group and $X$ is a scheme, and the categorical action is
induced from a usual group action $\nu: G\times X\to X$ via the
formula on objects $F\triangleright L = \nu_*(F\boxtimes L)$ where
$\boxtimes$ denotes the external tensor product of sheaves.

Leaving apart the difficult question on which monoidal categories $\mathcal G$
generalizing ${}_H\MM$ for Hopf algebra $H$, should qualify
as representing the 'noncommutative group schemes', 
we set the convention that $\mathcal G$ 
will be equipped with a distinguished action on
the base category (which in general does not need to be monoidal);
we consider only the actions respecting (via direct image
functors) the distinguished action on the base, and call such
actions {\bf (geometrically)} {\bf admissible} (or strictly
compatible with the distinguished action in the base). This
innocent condition is in fact very central to our approach! The
following result is now a geometric restatement of our simple
earlier result in \cite{Skoda:distr}:

\nxpoint {\bf Proposition.} {\it If $\mathcal C$ is monadic over
the base $\mathcal B$, that is $\mathcal C\cong\mathcal B^{\mathbf T}$,
where the monad $\mathbf T$ on the base is
the composition of the inverse and direct image functors of the
geometric morphism, then for a monoidal category $\mathcal G$
representing a symmetry object over the base $\mathcal B$, the
distributive laws between the distinguished action of $\,\mathcal
G$ on the base $\mathcal B$ and the monad $\mathbf T$ are in a
bijective correspondence with the geometrically admissible
monoidal actions of $\mathcal G$ on $\mathcal C$. }
\vskip .05in

This simple analogue (secretly, a generalization) of the classical
Beck's theorem (\cite{Borceux,Semtriples}) on the bijection
between the distributive laws of monads and lifts of one monad to
the Eilenberg--Moore category of another monad explains 
many appearances of ``entwining structures'' in noncommutative
fiber bundle theories.

\nxpoint \label{s:comodulealgebra}
Recall that the category of left modules over any
bialgebra is monoidal.
\vskip .06in

{\bf Definition.} {Let $B$ be a bialgebra. A {\bf right
$B$-comodule algebra} $E$ is a right $B$-comodule for which the
coaction $\rho:E\to E\otimes B$ is an algebra map. }
\vskip .1in

{\bf Proposition.} {\it Every right $B$-comodule algebra $E$
canonically induces a geometrically admissible action
of ${}_B\MM$ on ${}_E\MM$.
}
\vskip .07in

{\it Proof.} It is sufficient to write down a canonical
distributive law enabling the lifting of
the geometrically admissible ${}_B\MM$-action from
${}_\genfd\MM$ to ${}_E\MM$. Let $\lozenge:{}_\genfd\MM\times
{}_B\MM\to {}_\genfd\MM$ be the ``trivial'' tensor product action of
${}_B\MM$ on ${}_\genfd\MM$. The monad in question is of course
$E\otimes_\genfd$ and the distributive law has the components
$l_{E,M,Q}:E\otimes(M\lozenge Q)\to(E\otimes M)\otimes Q$ which
are given by the $\genfd$-linear extension of formulas $e\otimes
(m\otimes q)\mapsto \sum e_{(0)} \otimes m\otimes e_{(1)} q$,
where $e\in E$, $m\in M$, $q\in Q$, $\rho(e) =\sum e_{(0)}\otimes
e_{(1)}$ is the formula for coaction in the extended Sweedler
notation (\cite{Majid}). Easy calculations show
that $l_{E,M,Q}$ is indeed a distributive law.

\nodo{In fact, every mixed distributive laws between an algebra
and coalgebra, so-called entwining $\psi : E\otimes C\to C\times
E$ where $C$ is a coalgebra induces also an example of the
distributive law of our kind (but there are other examples as
well). Indeed, tensoring $\genfd$-modules with a coalgebra may be
viewed as a comonad. }

This way, the comodule algebras supply examples of
geometrically admissible actions; hence we view them as (a
class of) noncommutative $G$-spaces.

\nxpoint \label{s:modulealgebra}
{\bf Definition.} {A $B$-{\bf module algebra} is an algebra $A$
with a $B$-action $\triangleright$
satisfying the ``Leibniz rule'' $b\triangleright (a a')
= \sum (b_{(1)}\triangleright a) (b_{(2)}\triangleright a)$.
}
\vskip .05in

{\bf Proposition.} {\it Let $A$ be a left $B$-module algebra.
Then the monoidal category of right $B$-comodules
acts on ${}_A\MM$.}
\vskip .06in

 Again, for all $\genfd$-modules $M$ and
 $B$-comodules $Q$, one needs to write the components of the distributive law
 $l_{M,Q} : A\otimes(M\otimes Q)\to (A\otimes M)\otimes Q$.
 Indeed, the formula $a\otimes (m\otimes q) \mapsto \sum (q_{(1)}
 \triangleright_A a \otimes m)\otimes q_{(0)}$ does the job.

\section{Equivariant Sheaves in Noncommutative Geometry}

\nxpoint \label{p:Yonedatrick} 
Mumford (\cite{Mumford:GIT}, 1.3) defines equivariant
sheaves using an explicit cocycle condition. We start with a 
conceptually simple definition of an equivariant object 
in a fibered category (\cite{Vistoli}), which is easy to generalize.

Given a category $\mathcal C$ and an internal group $G$ in
$\mathcal C$, the Yoneda embedding induces a presheaf of groups
$h_G$ on $\mathcal C$. Given any presheaf of groups
$\hat{G}:\mathcal C^{\circ}\to\mathrm{\underline{Group}}$ over
$\mathcal C$, an action of $\hat{G}$ on an object $X$ in $\mathcal C$
is given by a natural transformation of functors
$\nu:\hat{G}\times h_X\to h_X$ such that for each object $U$ in
$\mathcal C$ the component $\nu_U :\hat{G}(U)\times\hom(U,X)\to\hom(U,X)$
is a group action of the group $G(U)$ on a set $\hom(U,X)$.
One obtains the category $G$-$\mathcal C$
of $G$-actions in $\mathcal C$.

Let now $\pi :\mathcal F\to\mathcal C$ be a fibered category and
$\nu$ an action of $\hat{G}$ on the  fixed object $X$. The
composition $\pi\circ\hat{G}$ is a presheaf of groups in $\mathcal
F$ so one can form the category of $\pi\circ\hat{G}$-actions in
$\mathcal F$, and this category clearly projects via naturally
induced projection $\pi'$ to the category $G-\mathcal C$ (in fact this
projection, for cartesian closed $\mathcal C$ and $\mathcal F$
is a fibered category as well). The fiber
$(\pi')^{-1}(\id_{(X,\nu)})$ is {\bf the category of equivariant
objects} in $\mathcal F$ over $(X,\nu)$.

\nxpoint Cartesian product with $G$ is in fact a monad in
$\mathcal C$ and $h_G\times h_X \cong h_{G\times X}$. Thus one in
fact induces a presheaf over $\mathcal C$ of monads in $\catset$
and for any presheaf of monads $\bf T$ one can do
similar trick as in \refpoint{p:Yonedatrick}
to define the Eilenberg--Moore fibered category
$\mathcal F^{\bf T}\to\mathcal C^{\bf T}$, which may be viewed as
the fibered category of equivariant objects. Unfortunately,
few monads in $\mathcal C$ can be replaced using Yoneda by
a presheaf of monads in $\catset$.

Similarly, for a functor of $V$-enriched categories $\pi :\mathcal
F\to\mathcal C$ one uses the enriched Yoneda lemma to define, for
any presheaf of (co)monads $\bf T$ in $V$ (e.g., tensoring with a
(co)group in the monoidal category $V$), the
$\bf T$-equivariant objects in $\mathcal F$ over a $\bf T$-module
$(X,\nu)$ in the base $\mathcal C$.

\nxpoint One can enrich categories to get 2-categories, and apply
the above mechanism. But Yoneda lemma should be better
replaced by the pseudo-Yoneda lemma for 2-functors in a pseudo
sense (contravariant version: 2-presheaves).
 C.~Hermida studied in~\cite{Hermida} a concept of
{\it 2-fibered 2-category}. For 1-fibrations one requires that
every arrow has a (strongly) Cartesian lift. Hermida requires 2
universal properties for liftings of 1-cells (1-Cartesian and
2-Cartesian 1-cells) and 2 universal properties for 2-cells in
order to call a 2-functor a 2-fibered 2-category. Given a weak (in
pseudo sense) 3-functor from the base 2-category
$\CC^{\mathrm{op}}$ to $2\cat$ one can perform a 2-categorical
analogue of Grothendieck construction to obtain a 2-fibered
2-category in the sense of Hermida. Vice versa, using a
2-categorical analogue of a cleavage one may represent 2-fibered
2-categories by pseudo-2-functors. Recall that if $Y$ is an
internal (monoidal) category in a 2-category $\CC$, then
$\hom_\CC(-,Y)$ is a usual internal (monoidal) category.

\nxsubpoint {\bf Definition.} {Let $G$ be an internal monoidal
category in the base 2-category $\CC$ of a 2-fibered category
$\pi:{\mathcal F}\to\CC$. Then $\hom(-,G)$ gives a representable
2-presheaf with values in the 2-category of (usual) monoidal
categories $\mathbf{MonCat}$. Let $\hat{G} = \hom(-,G)\circ\pi :
{\mathcal F}\to\cat$ be a 2-presheaf of monoidal categories over
$\mathcal F$. A {\bf $G$-equivariant (or 2-equivariant) object $\rho$
in $\mathcal F$} over an internal $G$-actegory $(X,\triangleright,\Psi,u)$ 
(where $\triangleright : G\times X\to X$ is a monodial action internal functor
in $\CC$, and $\Psi,u$ are the coherences) 
is a natural transformation of 2-functors
$\tilde\triangleright : \hat{G}\times \hom(-,\rho)\to
\hom(-,\rho)$ together with modifications $\tilde\Psi$,
$\tilde{u}$ of natural transformations of 2-functors
$\tilde\Psi:(-\otimes -)\tilde\triangleright -\tto
-\tilde\triangleright (-\tilde\triangleright -) :
\hat{G}\times\hat{G}\times\hom(-,\rho)\tto\hom(-,\rho)$,
$\tilde{u}:(-\otimes {\bf 1})\to -$ and satisfying the action
pentagon and unit triangle for every fixed argument $U\in \mathcal
F$, and such that there the action of $\hat{G}$ on $\hom(-\rho)$
is compatible via projection with the action of $G$ on $X$, i.e.,
$\pi(\tilde\triangleright) = \hom_\CC(-,\triangleright)$,
$\pi(\tilde\Psi) = \hom_\CC(-,\Psi)$,
$\pi(\tilde{u})=\hom_\CC(-,u)$. } \vskip .02in

The last part is just symbolic and may need clarification.
Given \mbox{$U\!\in\!\mathcal F$}, 
the $U$-component of the equality $\pi(\tilde\triangleright) =
\hom_\CC(-,\triangleright)$
means that for all $g\in \hat{G}(U)$, $z\in\hom(U,\rho)$,
$\pi(g\,\tilde\triangleright z)= \pi(g)\triangleright\pi(z)$ and
alike for the modifications $\tilde\Psi, \tilde{u}$.

This clear definition of 2-equivariant objects is new (we
anounced further results involving this definition at WAGP06
conference ``Gerbes, groupoids and QFT'' at ESI, Vienna, in May
2006). One can write explicit cocycle descriptions of
2-equivariant 2-objects in $\mathcal F$ in the style of Mumford.
In physics, actions of groups on gerbes, involve a special case of
2-(categorical) equivariance, which are discussed in
\cite{konradUrsYandl}. The detailed treatment and applications
will appear elsewhere.

\nxpoint \label{s:Hopfmodules} {\bf (Comonad for the relative Hopf
modules).} Let $B$ be a $\genfd$-bialgebra. To any right
$B$-comodule algebra $(E,\rho_E)$ (\refpoint{s:comodulealgebra}),
we associate an endofunctor $G : {}_E\MM\to{}_E\MM$ in the
category ${}_E\MM$ of left $E$-modules on objects $M$ in ${}_E\MM$
given by the formula $G: M\mapsto M\otimes B$, where the left
$E$-module structure on $M\otimes B$ is given by $e(m\otimes b):=
\rho_E(e)(m\otimes b)$, or in an extended Sweedler notation
(\cite{Majid}), $e(m\otimes b) = \sum e_{(0)} m\otimes e_{(1)}b$
($e\in E, m\in M, b\in B$). In calculations we will often write
just $e_{(0)}\otimes e_{(1)}$, omiting even the summation sign
$\sum$ in the Sweedler notation. The comultiplication
$\Delta=\Delta^B$ on $B$ induces the comultiplication $\delta
=\id\otimes\Delta:G\to GG$ on $G$ with counit
$\epsilon^G=\id\otimes\epsilon$ making
$\mathbf{G}=(G,\delta,\epsilon^G)$ a comonad (cf. the coring
picture in~\cite{BrzWis:corings}).

\nxsubpoint A left-right {\bf relative $(E,B)$-Hopf module} is a
triple $(N,\rho_N,\nu_N)$ such that $\nu_N:E\otimes N\to N$ is a
left $E$-action, $\rho_N:N\to N\otimes B$ is a right $B$-coaction
and
$\rho_N(\nu(e,n))=(\nu\otimes\mu_B)(\id\otimes\tau_{B,N}\otimes\id)
(\rho_E(e)\otimes\rho_N(n))$ for all $e\in E$, $n\in N$, where
$\tau_{B,N}:B\otimes N\to N\otimes B$ is the flip of tensor
factors. Maps of relative Hopf modules are morphisms of underlying
$\genfd$-modules, which are maps of $E$-modules and $B$-comodules.

\nxsubpoint {\bf Proposition.} {\it The category $({}_E\MM)_{\bf
G}$ of $\mathbf{G}$-comodules $($coalgebras$)$ is equivalent to the
category ${}_E\MM^B$ of left-right relative $(E,B)$-Hopf modules.
}
\vskip .013in

This is one of our basic observations 
in a collaboration with V. Lunts (2002),
and is independently observed and used for
similar purposes in coring theory about at the same time (and
generalizations for entwined modules and so on).

\nxpoint P. Deligne in \cite{Del:Hodge3} notes that the category
of $G$-equivariant sheaves naturally embeds into the category of
simplicial sheaves over the Borel construction considered as a
simplicial space. In Autumn 2002, we noticed with V.~Lunts
that a parallel construction exists for relative $(E,B)$-Hopf
modules.

\nxsubpoint Recall that any comonad ${\bf G}$ on a category
$\mathcal A$ induces an augmented {\it simplicial} endofunctor
${\mathbf G}_\bullet\stackrel\epsilon\to \Id_{\mathcal A}$ and
dually a monad induces an augmented cosimplicial endofunctor.
Starting with the comonad ${\bf G}$ from \refpoint{s:Hopfmodules}
on ${}_E\MM$, we can form the category of ${\bf G}$-comodules
$({}_E\MM)_{\mathbf G}\cong {}_E\MM^B$, and form a {\it monad}
${\mathbf T}_{\mathbf G}$ on $({}_E\MM)_{\mathbf G}$ obtained from
the adjoint pair of the forgetful and cofree functors
$U_G:({}_E\MM)_{\mathbf G}\leftrightarrow {}_E\MM:F_G$. Thus this
{\it monad} ${\mathbf T}_{\mathbf G}$ on $({}_E\MM)_{\mathbf G}$
induces a {\it cosimplicial} endofunctor ${\mathbf T}_{\mathbf
G}^\bullet$ on ${}_E\MM^B$. As $(E,\rho)$ is a monoid in
${}_E\MM^B$, ${\mathbf G}^\bullet (E,\rho) := U_{\mathbf G}
{\mathbf T}_{\mathbf G}^\bullet (E,\rho)$ is a cosimplicial
algebra in ${}_E\MM$, the {\bf coborel construction} on
$(E,\rho)$.

Let $f:{\bf n}\to{\bf m}$ be a morphism in the category of
nonempty finite ordinals (simplices) $\Delta$, and $G_f : G^n E\to
G^m E$ the corresponding map in $G^\bullet E$. The following idea
is due to V. Lunts:

\nxsubpoint {\bf Definition.} {A {\bf simplicial module}
$M_\bullet$ over the coborel construction ${\mathbf G}^\bullet
(E,\rho)$ is a sequence $(M_n)_{n = 0,1,2,\ldots}$ of
$\genfd$-modules, with the structure of a left $E_n$-module on
$M_n$, together with structure maps of left $E_n$-modules $\beta_f
: G_f^* M_m\to M_n$, for all morphisms $f :{\bf n}\to{\bf m}$ in
$\Delta$ (where $G_f^*$ is the extension of scalars along $G_f$),
and such that for all ${\bf n}\stackrel{f}\to{\bf
m}\stackrel{g}\to{\bf r}$ in $\Delta$ the cocycle condition
holds:}
$$\xymatrix{G_{g\circ f}^* M_r \ar@/_1pc/[rrr]_{\beta_{g\circ f}}\ar[r]^\cong &
G_f^* G_g^* M_r \ar[r]^{G_f^* (\beta_g)} &G_f^* M_m
\ar[r]^-{\beta_f} & M_n}$$

The morphisms of simplicial modules are ladders of maps $M_n\to
N_n$ of $E_n$-modules $n=0,1,2,\dots$ commuting with the structure
maps $\beta_f$. This way we get a category ${}_E \mathbf{Sim}^B$
of simplicial modules over $G^\bullet
(E,\rho)$.

\nxsubpoint \label{lunts:thmesimb}{\bf Theorem} (with V. Lunts).
{\it The category of relative Hopf modules ${}_E\MM^B$ is
equivalent to the full subcategory of ${}_E \mathbf{Sim}^B$ of
those objects for which all $\beta_f$ are isomorphisms. } \vskip
.03in

We found several interesting proofs of this theorem (our forthcoming 
paper with V. Lunts: {\it Hopf modules, {\tt Ext}-groups and descent}). 
One of the proofs is via an intermediate construction of independent interest:

\nxsubpoint Given a right $B$-comodule algebra $(E,\rho)$, denote
by $p = p_E : E\to E\otimes B$ the map $e\mapsto e\otimes 1_B$ and
by $p_{12}:E\otimes B\to E\otimes B\otimes B$ the map $e\otimes
b\mapsto e\otimes b\otimes 1$. A {\bf right $B$-coequivariant left
$E$-module} is a pair $(M,\theta)$ where $M$ is a left $E$-module
and $\theta : \rho^* M \rightarrow p^* M$ is an isomorphism of
left $E$-modules which satisfies the following Mumford-style
cocycle condition:
\begin{equation}\label{eq:cocycle}\xymatrix @C=3.1pc { (\id \otimes
\Delta)^* \rho^* M \ar[d]^\cong \ar[r]^{(\id\otimes\Delta)^*
\theta} & (\id \otimes \Delta)^* p^* M \ar[r]^-\cong &
 p_{12}^* p^* M \\
(\rho \otimes \id)^* \rho^* M\ar[r]^{(\rho\otimes\id)^*\theta} &
 (\rho \otimes \id)^* p^* M \ar[r]^-\cong & p_{12}^* \rho^* M
\ar[u]^{p_{12}^* \theta}}
\end{equation}
Notice that in our notation $f\mapsto f^*$ is a {\it covariant}
functor, and that the three canonical isomorphisms denoted by $\cong$
in the diagram are nontrivial when written in terms of tensor
products. A {\bf morphism of pairs} $f :
(M,\theta_M)\to(N,\theta_N)$ is a morphism of left $E$-modules
$f:M\to N$ such that $p^* f \circ \theta_M = \theta_N\circ
\rho^*f$. This way we obtain a category ${}_E\MM^{\mathrm{coeq}B}$
of right $B$-coequivariant left $E$-modules.

\nxsubpoint {\bf Theorem} (with V. Lunts). {\it There is a
canonical equivalence of categories ${}_E\MM^B\cong
{}_E\MM^{\mathrm{coeq}B}$.}

{\it Sketch of the proof.} If $f:E\to E\otimes B$ then the class
in $f^*M$ with representative $e\otimes m\otimes b$ in $(E\otimes
B)\otimes M$ will be denoted $[e\otimes m\otimes b]_{f^*M}$.

The equivalence of categories needs to produce $\theta$ from
$\rho_M$ and viceversa (the underlying $M\in{}_E\MM$ does not
change). Given coaction $\rho_M$, define $\theta : \rho^*_{E} M\to
p^* M$ by the $\genfd$-linear extension of the formula
\[ \theta([e \otimes b\otimes m]_{\rho^* M})
:= [\sum e \otimes bm_{(1)}\otimes m_{(0)}]_{p^* M}.\]

Given $\theta$, define $\rho_M : M \to M \otimes B$ by $\rho_M (m)
:= ({\rm nat}\,\circ \theta)[1 \otimes 1\otimes m]_{\rho^*_E M}$.
Map ${\rm nat} : p^* M \to M\otimes B$ is the isomorphism of
$\genfd$-linear spaces which is the composition in the bottom line
of the commutative diagram
\[\begin{array}{ccccc}
 E\otimes B \otimes M &
\stackrel{E\otimes \tau_{M,B}}\longrightarrow & E\otimes M \otimes
B& \stackrel{\nu\otimes B} \longrightarrow &
 M \otimes B \\
\downarrow & & \downarrow &&\downarrow \\
p^*M & \longrightarrow & \frac{E\otimes M \otimes B}{\langle e
\otimes m \otimes b - 1 \otimes e m \otimes b\rangle} &
\longrightarrow & M \otimes B
\end{array}\]
where $\nu$ is the action $E\otimes M \to M$ and $\tau_{M,B} :
B\otimes M \to M \otimes B$ is the flip of tensor factors and the
vertical lines are the natural projections.

Now all the verifications (the correspondences are well defined
and mutually inverse, the cocycle condition for the new theta, the
coaction axiom for new $\rho_M$) are just calculations with
classes $[e\otimes m\otimes b]_{f^* M}$ (and "longer" versions).
Finally, in the same terms, one checks that the map of left
$E$-modules $h:M\to N$ is a morphism $(M,\rho_M)\to(N,\rho_N)$ iff
it is morphism between the corresponding coequivariant modules
$(M,\theta_M)\to(N,\theta_N)$. This finishes the proof.\vskip
.03in

Theorem \refpoint{lunts:thmesimb} can now be proved along the
following lines: starting with $(M,\rho_M)$ we first form
$(M,\theta)$, then we set $M_0 = M$, $M_1 = p^* M$, $M_2 =
p_{12}^* p^* M$ and so on. Maps $p,p_{12},p_{123}$ etc. are the
$0$-th coface maps of the coborel construction, and the
corresponding structure morphisms can be taken identities. A
general structure map can be computed easily if we know the
structure maps corresponding just to cofaces and codegeneracies.
But those can be easily found from comparing domains $f^* M_m$
with $p_{12..}^* M_m$. For example, $\rho^* M_0 \to M_1$ is the
composition $\rho^* M_0 \stackrel\theta\longrightarrow p^*
M_0\stackrel{=}\to M_1$, that is simply $\theta$. After similar
simple formulas are proposed for all generators, one can check the
properties.

\nxsubpoint Above considerations propose the following more
general definition applicable to a large class of cases, including
relative Hopf modules and classical equivariant sheaves.

Let $\bf G$ be a (co)monad in a base category $\CC$ of a fibered
category  $\pi: \mathcal F\to\CC$ and $(E,\rho)$ a ${\bf
G}$-(co)module. The {\bf category of equivariant objects over $E$}
is the category of Cartesian functors from $\Delta^0$
(respectively $\Delta$) considered as a discrete fibered category
into $\mathcal F$, whose bottom part is the (co)bar construction
for $(E,\rho)$. In other words, it is the fiber (the category of
Cartesian sections) over the (co)bar construction considered as a
functor. For Hopf modules the base category is the category of
$\genfd$-algebras, the fiber over $A$ is ${}_A\MM$, the pullback
is the extension of scalars, and the cobar construction is our
coborel construction.

Similar constructions can be made for (co)lax actions of monoidal
categories, generalizing the (co)monad case.

\section{Compatible Localizations}

\nxpoint Here we consider flat localizations of rings
(e.g., Ore localizations), and also (additive) localization functors $Q^*$
(between Abelian categories) possesing a right adjoint $Q_*$ (we call them
continuous localization functors); equivalently, $Q_*$ is
a fully faithful functor having a left adjoint; or counit of the adjunction
is an isomorphism (equivalently, the multiplication of the
corresponding monad is an isomorphism, i.e., the monad is idempotent)
(\cite{GZ,skoda:ncloclec}).


The following concept has been introduced in my thesis (the thesis
results are published in
\cite{Skoda:ban,Skoda:Ore,Skoda:coh-states}).

\nxpoint {\bf (Compatibility of coactions and localizations).}
(\cite{Skoda:ban}) Given a bialgebra $B$ and a (say, right)
$B$-comodule algebra $(E,\rho)$, an Ore localization of rings
$\iota : E\to S^{-1}E$ is $\rho$-{\bf compatible} if there exists
an (automatically unique) coaction $\rho_S : S^{-1}E\to
S^{-1}E\otimes B$ making $S^{-1}E$ a $B$-comodule algebra, such
that the diagram
$$\xymatrix{
E\ar[r]^\rho\ar[d]^{\iota_S}& E\otimes B\ar[d]^{\iota_S\otimes B}
\\
S^{-1}E \ar[r]^{\rho_S}&S^{-1}E\otimes B }$$ commutes; $\rho_S$ is
then called the {\bf localized coaction}. We call the $\rho_S$-coinvariants in $S^{-1}E$ {\bf localized coinvariants}. Even for
compatible localizations, $\iota_S$ restricted to the subring
$E^{\mathrm{co}B}\subset E$ is typically not underlying the ring
localization $U^{-1} E^{\mathrm{co}B}$ with respect to any Ore
subset $U$ in $E^{\mathrm{co}B}$.

\nxpoint \label{s:compOreCateg}{\bf Theorem.} {\it Let $B$ be a
$\genfd$-bialgebra, $(E,\rho)$ a $B$-comodule algebra,
$\mathbf{G}$ the corresponding comonad on the category of left
$E$-modules, which is described in~\refpoint{s:Hopfmodules}, and
let $\iota : E\to E_\mu$ be a perfect $($e.g., Ore$)$ localization of
rings, which happens to be $\rho$-compatible. The $\genfd$-linear
maps
$$
l_M : E_\mu\otimes_E(M\otimes B)\!\rightarrow \!(E_\mu\otimes_E
M)\otimes B, \;\; e\otimes(m\otimes b)\!\mapsto\! \sum(e_{(0)}\otimes
m)\otimes e_{(1)} b,
$$
where $M$ runs through left $E$-modules, are well-defined morphisms
of left $E$-modules and together they form a mixed distributive
law $l : Q_* Q^* G\tto G Q_* Q^*$ between the localization monad
$Q_* Q^*$ and the comonad $\mathbf G$ on ${}_E\MM$. }

\begin{proof}
Clearly, the ($\genfd$-linear extension of the) formula
$l':e\otimes(m\otimes b)\mapsto \sum (e_{(0)}\otimes m)\otimes
e_{(1)} b$ gives a well-defined $\genfd$-linear map
$$
E_\mu\otimes_\genfd M\otimes_\genfd B\to E_\mu\otimes_\genfd
M\otimes_\genfd B.
$$
To show that $l'$ factors to a well-defined map
$l : E_\mu\otimes_E (M\otimes_\genfd B)\to (E_\mu\otimes_E
M)\otimes_\genfd B$
(where the $E$-module structure on $M\otimes B$ is from
\refpoint{s:Hopfmodules}), we need to show that if
$e'\in E_\mu$ and $e\in E$, then the  map $l'$ sends
$r = e' e\otimes m\otimes b$ and
$s = \sum e'\otimes e_{(0)}m\otimes e_{(1)}b$
to the elements in the same class in $(E_\mu\otimes_E M)\otimes B$.
This is easy:
$l'(r) = \sum ((e'e)_{(0)}\otimes m)\otimes
(e'e)_{(1)} b = \sum (e'_{(0)}e_{(0)}\otimes m)\otimes
e'_{(1)}e_{(1)} b$ and
$l'(s) = \sum (e'_{(0)} \otimes e_{(0)}m)\otimes e'_{(1)}e_{(1)}b$,
which clearly becomes the same class when projected to
$(E_\mu\otimes_E M)\otimes B$
(move $e_{(0)}$ along $\otimes_E$).
Thus $l_M$ is a well-defined $\genfd$-linear map.
It is easy to see that $l_M$ is also $E$-linear (and even $E_\mu$-linear).

To check the first pentagon
$$\xymatrix{
Q_* Q^* G \ar[d]\ar[rr] && G Q_* Q^*\ar[d]\\
Q_* Q^* GG\ar[r]& GQ_* Q^* G\ar[r]& GG Q_* Q^*
}$$
for some $M$ in ${}_E\MM$,
we just directly calculate the two paths starting from
a generic homogeneous element
$e\otimes (m\otimes b)\in E_\mu\otimes_E (M\otimes B) = Q_* Q^* GM$.
{\footnotesize
$$
\xymatrix{
e\otimes(m\otimes b) \ar@{|->}[d]\ar@{|->}[rr] && (e_{(0)}\otimes
m)\otimes e_{(1)}b_{(1)}\ar@{|->}[d]\\
e\otimes ((m\otimes b_{(1)})\otimes b_{(2)})\ar@{|->}[r]&
(e_{(0)}\otimes (m\otimes b_{(1)}))\otimes e_{(1)} b_{(2)}
\ar@{|->}[r]& e_{(0)}\otimes m\otimes e_{(1)} b_{(1)}\otimes e_{(2)} b_{(2)}
}$$}

The second pentagon
$$
\xymatrix{
Q_* Q^* Q_* Q^*\ar[r]^{Q_* Q^* l}\ar[d]^{Q_*\epsilon Q^*G}
&Q_* Q^* G Q_* Q^*\ar[r]^{l Q_* Q^* }&
GQ_* Q^* Q_* Q^*\ar[d]_{GQ_*\epsilon Q^*}\\
Q_* Q^*G\ar[rr]^{l}&&GQ_*Q^*
}$$
is calculated similarly on $f\otimes(e\otimes(m\otimes b))\in
E_\mu\otimes_E(E_\mu\otimes_E(M\otimes B))=Q_* Q^* Q_* Q^* GM$:
the upper and right arrows compose
$$\begin{array}{l}
f\otimes(e\otimes(m\otimes b))          \longmapsto
f\otimes ((e'_{(0)}\otimes m)\otimes e_{(1)} b) \\
\longmapsto (f_{(0)}\otimes (e_{(0)}\otimes m))\otimes f_{(1)}(e_{(1)}b)
\longmapsto (f_{(0)} e_{(0)}\otimes m)\otimes f_{(1)}e_{(1)} b
\end{array}$$
while the left-below path gives
$$\begin{array}{l}
f\otimes(e\otimes(m\otimes b))\longmapsto fe\otimes(m\otimes b) 
\\
\longmapsto ((fe)_{(0)}\otimes m)\otimes (fe)_{(1)} b
= (f_{(0)} e_{(0)}\otimes m)\otimes f_{(1)}e_{(1)} b
\end{array}$$
The component of the
$\bf G$-counit triangle $(\epsilon^G Q_* Q^*)\circ l = Q_* Q^*\epsilon^G$
at object $M$ in ${}_E\MM$, can be computed at any tensor monomial
$e\otimes m\otimes b\in E_\mu\otimes_E (M\otimes B) =  Q_* Q^* G M$.
The left-hand side computes as
$e\otimes m\otimes b\mapsto e_{(0)}\otimes m\otimes e_{(1)} b\mapsto e_{(0)}\otimes m\epsilon(e_{(1)}b) = e\otimes m\epsilon(b)$
and the right-hand side gives the latter result immediately.

Finally, the unit triangle
$$\xymatrix @R=1pc{
&GM\ar[rd]^{\eta_{GM}}\ar[ld]_{G(\eta_M)}&\\
E_\mu\otimes_E GM\ar[rr]^{l_M}&& G(E_\mu\otimes_E M) }$$ is almost
trivial to check: $m\otimes b\mapsto 1\otimes (m\otimes b)\mapsto
(1\otimes m)\otimes b$ and the direct arrow $GM\to
G(E_\mu\otimes_E M)$ is $m\otimes b\mapsto(1\otimes m)\otimes b$.
\end{proof}

\nxpoint {\bf Proposition.} {\it If $B$ is a Hopf algebra with
antipode $S:B\to B$, then the formula $l_M^{-1} : (e\otimes
m)\otimes b\mapsto e_{(0)}\otimes (m\otimes S(e_{(1)})b)$ defines
a $\genfd$-linear map $l_M^{-1}: (E_\mu\otimes_E M)\otimes B\to
E_\mu\otimes_E (M\otimes B)$ which is inverse to $l_M$ and is a
map of $E_\mu$-modules. }

\begin{proof}
The inverse property is obvious: $e\otimes (m\otimes b)
\stackrel{l_M}\longmapsto (e_{(0)}\otimes m)\otimes e_{(1)}b
\stackrel{l_M^{-1}}\longmapsto
(e_{(0)(0)}\otimes m)\otimes S(e_{(0)(1)})e_{(1)}b=e\otimes m\otimes b$.
The map $l_M^{-1}$ is indeed a homomorphism of $E_\mu$-modules, because
for $f\in E_\mu$ it sends $f((e\otimes m)\otimes b)
= (f_{(0)} e\otimes m)\otimes f_{(1)}b$,
by definition, to
$(f_{(0)} e)_{(0)}\otimes (m\otimes S(f_{(1)} e_{(0)})f_{(1)}b)
= f_{(0)}e_{(0)}\otimes (m\otimes Se_{(1)} Sf_{(1)}\cdot f_{(2)}b) =
fe_{(0)} \otimes (m\otimes e_{(1)}b)
= f(l_M^{-1}(e\otimes (m\otimes b)))$, as required.
\end{proof}

\nxpoint \label{s:lemmaLiftingGeneral} {\bf Lemma.} {\it Any map
of comonads $($possibly in different categories$)$ induces a functor
between their respective categories of comodules $($coalgebras$)$.
}

This follows from the dual statement for monads; the latter is a part
of a stronger fact that taking the Eilenberg--Moore category
extends to a (strict) 2-functor from the category of monads
to the category of categories, see \cite{Street:formal}.

\nxpoint \label{s:mapComonads}
{\bf Theorem.} {\it Given any continuous localization
functor $Q^*:\AA\to\AA_\mu$ and a comonad $\bf G$ together with
any mixed distributive law $l:Q_* Q^* G\tto G Q_* Q^*$,

{\rm (i)} $G_\mu = Q^* G Q_*$ underlies a comonad ${\bf
G}_\mu=(G_\mu,\delta^\mu, \epsilon^{G_\mu})$ on $\AA_\mu$ with
comultiplication $\delta^\mu$ given by the composite
$$\xymatrix @C3pc{
Q^* G Q_* \ar[r]^-{Q^* \delta^G Q_*}& Q^* GG Q_* \ar[r]^-{Q^*
G\eta G Q_*}& Q^* G Q_* Q^* G Q_* }$$ and whose counit
$\epsilon^{G_\mu}$ is the composite
$$\xymatrix @C4pc{
Q^*GQ_* \ar[r]^{Q^*\epsilon^G
Q_*}&Q^*Q_*\ar[r]^{\epsilon}&\Id_{\AA_\mu} }$$
$($where the
right-hand arrow $\epsilon$ is the counit of the adjunction
$Q_*\dashv Q^*)$.

{\rm (ii)} the composite
$$\xymatrix{
Q^* GM\ar[rr]^-{Q^*(\eta_{GM})}&& Q^* Q_* Q^* GM
\ar[rr]^-{Q^*(l_M)}&& G_\mu Q^* M }$$ defines a component of a
natural transformation $\alpha = \alpha_l: Q^* G\tto G_\mu Q^*$
for which the mixed pentagon diagram of transformations
\allowdisplaybreaks
$$\xymatrix{
Q^* G\ar[rr]^{\alpha}\ar[d]_{Q^*\delta^G}&& G_\mu Q^*\ar[d]^{\delta^\mu Q^*}\\
Q^* GG \ar[r]^{\alpha G}&G_\mu Q^* G\ar[r]^{G\alpha}& G_\mu G_\mu
Q^* }$$ commutes and $(\epsilon^{G_\mu} Q^*)\circ\alpha = Q^*
\epsilon^G$. In other words, $(Q^*,\alpha_l): (\AA, {\bf G})\to
(\AA_\mu, {\bf G}_\mu)$ is $($up to orientation convention which
depends on an author$)$ a map of comonads
$($\cite{Skoda:eqdis,Street:formal}$)$. } \vskip .1in

Part (i) is a standard general nonsense on distributive laws once
the continuous localization is replaced by the correspoding
idempotent monad. The proof of (ii) is easy.

\nxpoint {\bf Proposition.} {\it
Suppose $\iota:E\to E_\mu$ is the $\rho$-compatible
localization of a $B$-comodule algebra $E$, $\AA = {}_E\MM$,
$\AA_\mu={}_{E_\mu}\MM$, $\bf G$, $\bf \tilde{G}_\mu$
are the ``comonads for Hopf modules'' as in \refpoint{s:Hopfmodules} and
$G_\mu$, $l$ are constructed as in \refpoint{s:compOreCateg}.
The comonad $\bf G_\mu$ is isomorphic to the
comonad $\tilde{\bf G}_\mu$. Moreover $Q_* \tilde{G}_\mu$. $Q_* G_\mu$ and
$G Q_*$ are isomorphic endofunctors in ${}_{E_\mu}\MM$.
}

\begin{proof}
As  $\genfd$-vector spaces, clearly both
$Q_* \tilde{G}_\mu N$ and $G Q_* N$ look for $N\in {}_{E_\mu}\MM$ like
$N\otimes B$.
The $E$-module structure on $Q_* G_\mu N$ is restriction of
the $E_\mu$-module structure given by
$f(n\otimes b) =\rho_{E_\mu}(f)(n\otimes b)$ for $f\in E_\mu$,
that is $\iota(e)(n\otimes b)= \rho_{E_\mu}(\iota(e))(n\otimes b)$,
while the $E$-module structure on $GQ_* N$ is given by
$e(n\otimes b) = ((\iota\otimes\id_B)\rho_E(e))(n\otimes b)$.
By the $\rho$-compatibility of localization $\iota$ the two answers
agree, i.e., $Q_* \tilde{G}_\mu= G Q_*$.

$\epsilon : Q^* Q_*\tto \Id$ is an isomorphism. Hence
$G_\mu = Q^* G Q_* = Q^* Q_* \tilde{G}_\mu \cong \tilde{G}_\mu$.

One should further check that the comultiplications of $\bf G_\mu$ and
$\bf\tilde{G}_\mu$ agree, that is the external square in the diagram
$$\xymatrix{
E_\mu\otimes_E({}_E N\otimes B))\ar[d]^= \ar[rr]^{Q^* \delta^G Q_*}
&& E_\mu\otimes_E (({}_E N\otimes B)\otimes B)\ar[d]^= \\
Q^* Q_* ({}_{E_\mu}N\otimes B)
\ar[rr]^{Q^* Q_* \delta^{\tilde{G}_\mu}}\ar[d]^\epsilon
&& Q^* Q_*(({}_{E_\mu} N\otimes B)\otimes B)\ar[d]^\epsilon\\
{}_{E_\mu}N\otimes B\ar[rr]^{\delta^{\tilde{G}_\mu}}
&& ({}_{E_\mu} N\otimes B)\otimes B
}$$
commutes, where the upper vertical arrows are induced by isomorphisms
$Q_* \tilde{G}_\mu\!\cong G Q_*$. The lower square commutes by naturality
of $\epsilon$ and the upper square commutes as vertical
arrows are identities and horizontal arrows are both
$\id\otimes\Delta$ at the level of vector spaces.
\end{proof}

\nxpoint \label{s:complocLifting} {\bf Theorem.} {\it Under
the assumptions in \refpoint{s:compOreCateg}, there is a unique
induced continuous localization functor $Q^{B*}: {}_E\MM^B\to
{}_{E_\mu}\MM^B$ between the categories of relative Hopf modules
such that $U_\mu Q^{B*} = Q^* U$ where $U$ and $U_\mu$ are the
forgetful functors from the category of relative Hopf modules to
the categories of usual modules over $E$ and $E_\mu$, respectively.
}
\vskip .02in

\begin{proof}
This follows from \refpoint{s:compOreCateg} by
Lemma \refpoint{s:lemmaLiftingGeneral},
after applying the equivalences of categories
${}_E\MM^B\cong ({}_E\MM)_{\bf G}$ and ${}_E\MM^B\cong
({}_{E_\mu}\MM)_{{\bf G}_\mu}$.
\end{proof}

In~\cite{Skoda:lcomp}, we prove
an interesting generalization of Theorem \refpoint{s:complocLifting}
for $\rho$-compatible localizations to {\it localization-compatible
pairs of entwining structures} introduced therein.

\nxpoint The main reason why compatible localizations are needed
in noncommutative geometry is that they are the analogues of {\it
$G$-invariant} open sets (unions of $G$-orbits) in commutative
geometry, where $G$ is a group.

\nxpoint Lunts and Rosenberg studied (\cite{LuntsRosSel,
LuntsRosMP}) the rings of differential operators for
noncommutative rings, generalizing the commutative Grothendieck's
definition in a nontrivial way (another approach yielding the same
definition is in \cite{Maszczyk}). Their purpose was to generalize
the Beilinson--Bernstein localization theorem in representation
theory to quantum groups. The basic property of the differential
operators is that they extend to exact localizations. Beilinson and Bernstein
abstracted this to the compatibility of localization functors and
monads; and prove that it is satisfied for their basic object in
\cite{LuntsRosMP}, {\it differential monads}. The $D$-affinity of
Beilinson has its abstract and simple generalization in their
general context.

The compatibility between a comonad $\bf G$ and a continuous
localization functor $Q^*:\AA\to\BB$ of Lunts and Rosenberg
is an isomorphism of functors of the form $Q^* G\cong G' Q^*$
where $G'$ is some endofunctor in $\BB$. This looks like
our distributive law $Q^* G\tto G_\mu Q^*$. There are two differences:
our map is a distributive law (satisfies two pentagons and two triangles,
in which sense our definition is stronger), and they require an isomorphism
while we have only a natural transformation
(here their definition is stronger; in our main examples,
induced by comodule algebras, we get invertibility
for the Hopf algebra case, while not for bialgebras).

\nxpoint The notion of the compatibility of (co)monads and localization
functors can easily be extended to the compatibility of
actions of monoidal groups and localizations.
For simplicity, we leave out this generalization
(and present it in \cite{Skoda:lcomp}).
We implicitly (verbally in a definition) use it in the next section though.

\section{Principal Bundles and Quotient Schemes}

\nxpoint In commutative algebraic geometry, there is a notion of
the descent along torsors: given a group scheme $G$, the category
$\qcoh^G(Y)$ of $G$-equivariant quasicoherent sheaves on the total
space of a $G$-torsor $Y$ (in fpqc topology) over a scheme $X$ is
equivalent to the category $\qcoh(X)$ of usual quasicoherent
sheaves over $X$. More generally, take any site $\CC$, a group
object $G$ in $\CC$, $G$-torsor $Y$ over $X$, and any stack of
categories over $\CC$ (replacing the stack $\mathcal F$ of
categories of quasicoherent sheaves over the site of schemes in
fpqc topology); then $G$-equivariant fiber ${\mathcal F}^G_Y$ over
$Y$ is equivalent to the usual fiber ${\mathcal F}_X$ over $X$
(\cite{Vistoli}). The analogue holds for $(E,B)$-Hopf modules: the
theorem of Schneider (\cite{Schneider}) states that, given a
faithfully flat Hopf-Galois extension $U\hookrightarrow E$, the
category of relative $(E,B)$-Hopf modules is equivalent to the
category of left modules over $U$; this theorem has many
generalizations for entwining modules and so on. A {\bf
Hopf-Galois extension} is the inclusion $U\hookrightarrow E$ of an
algebra $U=E^{\mathrm{co}B}$ of $B$-coinvariants in a $B$-comodule
algebra $E$ into $E$, such that the canonical map $E\otimes_U E\to
E\otimes B$, $e\otimes e' \mapsto \sum e_{(0)}e'\otimes e_{(1)}$
is an isomorphism.

\nxpoint However, already in the commutative geometry we know that
it is a rare case that there are sufficiently many coinvariants to
reconstruct the quotient of an affine torsor under a group action
(the spectrum of the algebra of coinvariants is sometimes called
the {\it affine quotient}). Thus I suggest below a globalization of
Hopf-Galois extensions to the noncommutative schemes of Rosenberg
(\cite{Ros:NcSch,Ros:lecs}).\vskip .02in

A flat {\bf cover} $\{F_\mu : {\mathcal A}\to {\mathcal
B}_\lambda\}_{\lambda\in\Lambda}$ is a {\it conservative} family
of flat functors (i.e., a morphism $a:A\to A'$ in ${\mathcal A}$ is
invertible iff $F_\mu(a)$ is invertible for every
$\lambda\in\Lambda$). Recall the terminology of
\refpoint{s:pseudo}.

\nxpoint \label{s:ncscheme} A {\bf quasicompact relative
noncommutative scheme} $({\mathcal A},{\mathcal O})$ over a category
${\mathcal V}$ as an abelian category ${\mathcal A}$ with a distinguished
object ${\mathcal O}$, finite biflat affine cover by localizations
$Q^*_\lambda : {\mathcal A}\to {\mathcal B}_\lambda$, with a continuous
morphism $g$ from ${\mathcal A}$ to ${\mathcal V}$ (think of it as $X\to
{\rm Spec}\,\genfd$) such that each $g_* \circ Q_{\lambda *} :
{\mathcal B}_\lambda\to{\mathcal V}$ is affine. If ${\mathcal V} =
\genfd-$${\rm Mod}$ then ${\mathcal O} = g^*(\genfd)$ where $\genfd$
is viewed as an object in $\genfd-{\rm Mod}$. For more details and
examples see~\cite{Ros:NcSch}.

\nxpoint Let me propose a globalization
of the notion of Hopf-Galois extension:\vskip .02in

{\bf Definition.} {Given a Hopf algebra $B$, a noncommutative
scheme $\AA$ over ${}_\genfd\MM$ is a {\bf noncommutative
$B$-torsor} (over appropriate quotient) if there is a
geometrically admissible action $\lozenge:{}_B\MM\times\AA\to\AA$,
and an affine flat cover of $\AA$ by localizations
$Q^*_\mu:\AA\to\AA_\mu$, with $\AA_\mu\cong{}_{E^\mu}\MM$; where
each $Q^*_\mu$ is compatible with the comonad ${\bf G}$ induced by
$B$ (viewed as a comonoid in ${}_B\MM$) and action $\lozenge$; and
for each $\mu$, the induced comonad ${\bf G}_\mu$ in $\AA_\mu$, is
induced by a $B$-comodule algebra on $E^\mu$, and
$(E^\mu)^{\mathrm{co}B}\hookrightarrow E^\mu$ is a faithfully flat
Hopf-Galois extension. } \vskip .06in

Given a noncommutative $B$-torsor, algebras of localized
coinvariants $(E^\mu)^{\mathrm{co}B}$ are local coordinate
algebras of a cover of the quotient space represented by the
Eilenberg-Moore category $\AA^\lozenge$ of the monoidal action
$\lozenge$ (${}_B\MM$-equivariant sheaves on the space represented
by $\AA$). In other words, the definition explores a local version
of the Hopf-Galois condition, which allows the local
coordinatization of the ``quotient stack'' $\AA^\lozenge$. My main
example (\cite{Skoda:ban}) is the quantum fibration $SL_q(n)\to
SL_q(n)/B_q(n)$, (or the $GL_q(n)$ version $GL_q(n)\to
GL_q(n)/\tilde{B}_q(n)$, which has multiparametric deformation
generalizations), where the noncommutative spaces
$SL_q(n),GL_q(n)$ are represented by the quantum linear groups
${\mathcal O}(SL_q(n)),{\mathcal O}(GL_q(n))$, and the quotient
stack $SL_q(n)/B_q(n)$ is constructed as a noncommutative scheme
in my earlier work \cite{Skoda:ban}. This ``quantum fibration'' is
a ${\mathcal O}(B_q(n))$-torsor with local trivialization given by
$n!$ Hopf-Galois extensions in an interesting cover by $n!$ Ore
localizations (defined in terms of quantum minors)
$S_w^{-1}{\mathcal O}(SL_q(n))$ ($w$ in the Weyl group), which is an
analogue of the standard cover of $SL_n$ by shifts of the main
Bruhat cell. Our calculational proof shows existence of such a
cover by formulas using the quantum Gauss decomposition,
calculations with quantum minors, an explicit check of a
nontrivial Ore property (\cite{Skoda:Ore}), a difficult check of
the compatibility of coactions and localizations and, finally, the
local Hopf-Galois property which is in this case implied by the
stronger ``local triviality'', meaning that there is a smash
product decomposition $S_w^{-1}{\mathcal O}(SL_q(n)) \cong
U_w\sharp B_q$ (induced by quantum Gauss decomposition) as a right
$B_q$-comodule algebra. Embedding of coinvariants
$U_w\hookrightarrow U_w\sharp B_q$ into the smash product is a
simple case of the Hopf-Galois extension. This picture was applied
in (\cite{Skoda:coh-states}) to develop a geometric theory of
Perelomov coherent states for quantum groups with nontrivial
resolution of unity formula obtained utilizing line bundles over
the quantum coset space $SL_q(n)/B_q(n)$.

\section*{Acknowledgements}
The influence of discussions with V. Lunts,  A. L. Rosenberg and
M. Jibladze are widely present throughout my work on the subject
of the article. This article was written at the Rudjer
Bo\v{s}kovi\'c Institute in Zagreb and Max Planck Institute for
Mathematics in Bonn whom I thank for support. The travel to Bonn
was supported by DAAD Germany/MSES Croatia bilateral project.

\

\

\centerline{(Received 28.11.2008; revised 20.01.2009)}

\

Author's address:
\medskip

Division of theoretical physics

Institute Rudjer Bo\v{s}kovi\'c

P.O.Box 180, Zagreb HR-10002

Croatia

E-mail: zskoda@irb.hr

\begin{thebibliography}{99}
\bibitem{Semtriples}
{\sc H.~Appelgate, M.~Barr, J.~Beck, F.~W. Lawvere, F.~E. J. Linton,
E,~Manes, M.~Tierney,} and {\sc F.~Ulmer},
Seminar on triples and categorical homology theory (ETH 1966/1967). Edited by B. Eckmann. \emph{Lecture Notes in Mathematics,} No. 80. \emph{Springer-Verlag, Berlin--New York,} 1969.

\smallskip \bibitem{Bakovic}
{\sc I.~Bakovi\'c}, {Bigroupoid 2-torsors}. \emph{PhD thesis, LMU Munich,} 2008.
\newline
{\tt http://edoc.ub.uni-muenchen.de/9209}

\smallskip \bibitem{Berstein:cogroups}
{\sc I. Berstein}, {On cogroups in the category of graded
algebras}. \emph{Trans. Amer. Math. Soc.} {\bf 115}(1965), 257--269.

\smallskip \bibitem{Bohm:ent}
{\sc G.~B\"ohm}, Internal bialgebroids, entwining structures and corings. \emph{Algebraic structures and their representations,} 207--226, \emph{Contemp. Math.,} 376, \emph{Amer. Math. Soc., Providence, RI,} 2005, {\tt
math.QA/0311244.}

\smallskip \bibitem{Borceux}
{\sc F.~Borceux}, {Handbook of categorical algebra}, 3 vols.
\emph{Encyclopedia of Mathematics and its Applications,} 50, 51, 52. \emph{Cambridge University Press, Cambridge,} 1994. 

\smallskip \bibitem{BrzMaj:quantumgauge}
{\sc T.~Brzezi\'nski} and {\sc S.~Majid}, Quantum group gauge theory on quantum spaces. \emph{Comm. Math. Phys.} \textbf{157}(1993), No. 3, 591--638.

\smallskip \bibitem{BrzMajid:ent}
{\sc T.~Brzezi\'nski} and {\sc S.~Majid}, Coalgebra bundles. \emph{Comm. Math. Phys.} \textbf{191}(1998), No. 2, 467--492.

\smallskip \bibitem{BrzWis:corings}
{\sc T.~Brzezi\'nski, R.~Wisbauer}, Corings and comodules. \emph{London Mathematical Society Lecture Note Series,} 309. \emph{Cambridge University Press, Cambridge,} 2003.

\smallskip \bibitem{Del:Hodge3}
{\sc P.~Deligne}, Th\'{e}orie de Hodge. III. \emph{Inst. Hautes \'{E}tudes Sci. Publ. Math.} No. 44 (1974), 5--77.

\smallskip \bibitem{Del:tan}
{\sc P.~Deligne}, Cat\'{e}gories tannakiennes.
\emph{The Grothendieck Festschrift,} Vol. II, 111--195, \emph{Progr. Math.,} 87, \emph{Birkh\"{a}user Boston, Boston, MA,} 1990.

\smallskip \bibitem{Dem-Gab-eng}
{\sc M.~Demazure} and {\sc P.~Gabriel},
{\em Groupes alg\'{e}briques. Tome I.} \emph{Masson $\&$ Cie, Editeur, Paris; North-Holland Publishing Co., Amsterdam,} 1970; English transl.:
Introduction to algebraic geometry and algebraic groups.  \emph{North-Holland Mathematics Studies,} 39. \emph{North-Holland Publishing Co., Amsterdam--New York,} 1980.


\smallskip \bibitem{Drinf:ICM}
{\sc V.~G. Drinfel'd}, Quantum groups. \emph{Proceedings of the International Congress of Mathematicians,} Vol. 1, 2 (\emph{Berkeley, Calif.,} 1986), 798--820, \emph{Amer. Math. Soc., Providence, RI,} 1987.

\smallskip \bibitem{Fresse:cogroups} {\sc B. Fresse}, Cogroups in algebras over an operad are free algebras. \emph{Comment. Math. Helv.} \textbf{73}(1998), No. 4, 637--676.

\smallskip \bibitem{GZ}
{\sc P. Gabriel, M. Zisman},
Calculus of fractions and homotopy theory. \emph{Ergebnisse der Mathematik und ihrer Grenzgebiete,} Band 35. \emph{Springer-Verlag, New York,} 1967.

\smallskip \bibitem{Hermida}
{\sc C. Hermida}, Descent on 2-fibrations and strongly 2-regular 2-categories. \emph{Appl. Categ. Structures} \textbf{12}(2004), No. 5-6, 427--459.

\smallskip \bibitem{KJ:actions}
{\sc G.~Janelidze} and {\sc G. M. Kelly}, A note on actions of a monoidal category. \emph{CT2000 Conference} (\emph{Como}). \emph{Theory Appl. Categ.} \textbf{9}(2001/02), 61--91 (electronic).

\smallskip \bibitem{KontsRos}
{\sc M.~Kontsevich} and {\sc A.~L. Rosenberg}, Noncommutative smooth spaces. \emph{The Gelfand Mathematical Seminars,} 1996--1999, 85--108, \emph{Gelfand Math. Sem., Birkh\"{a}user Boston, Boston, MA,} 2000,
{\tt arXiv:math.AG/9812158}.

\smallskip \bibitem{Lack:laxlimits}
{\sc S. Lack}, Limits for lax morphisms. \emph{Appl. Categ. Structures} \textbf{13}(2005), No. 3, 189--203.

\smallskip \bibitem{LuntsRosMP}
{\sc V.~A. Lunts, A.~L. Rosenberg}, {Differential calculus
in noncommutative algebraic geometry}. Max Planck Institute Bonn preprints:

{\em I. D-calculus
on noncommutative rings}, MPI 96-53;

{\em II.  D-calculus
in the braided case. The localization of quantized enveloping
algebras}, MPI 96-76, Bonn 1996. 

\smallskip \bibitem{LuntsRosSel}
{\sc V.~A. Lunts, A.~L. Rosenberg},
Differential operators on noncommutative rings. \emph{Selecta Math. (N.S.)} \textbf{3}(1997), No. 3, 335--359;
sequel: Localization for quantum groups. \emph{Selecta Math. (N.S.)} \textbf{5}(1999), No. 1, 123--159.

\smallskip \bibitem{MacLane}
{\sc S.~Mac Lane}, Categories for the working mathematician. \emph{Graduate Texts in Mathematics,} Vol. 5. \emph{Springer-Verlag, New York--Berlin,} 1971.

\smallskip \bibitem{Majid}
{\sc S. Majid}, Foundations of quantum group theory. \emph{Cambridge University Press, Cambridge,} 1995.

\smallskip \bibitem{Maszczyk}
{\sc T.~Maszczyk}, {Noncommutative geometry through monoidal categories I}. \\{\tt arXiv:0707.1542}.

\smallskip \bibitem{Mumford:GIT}
{\sc D. Mumford},
Geometric invariant theory. \emph{Ergebnisse der Mathematik und ihrer Grenzgebiete, Neue Folge,} Band 34. \emph{Springer-Verlag, Berlin--New York}, 1965.
Second edition with {\sc J. Fogarty}. \emph{Ergebnisse der Mathematik und ihrer Grenzgebiete} [\emph{Results in Mathematics and Related Areas}], 34. \emph{Springer-Verlag, Berlin,} 1982.

\smallskip \bibitem{ParshallWang}
{\sc B. Parshall} and {\sc J.Wang}, Quantum linear groups. \emph{Mem. Amer. Math. Soc.} \textbf{89}(1991), No. 439, vi+157 pp.

\smallskip \bibitem{Ros:lecs}
{\sc A.~L. Rosenberg}, {Topics in noncommutative algebraic
geometry, homological algebra and K-theory}. \emph{Preprint} MPIM2008-57,
at {\tt http://www.mpim-bonn.mpg.de}.

\smallskip \bibitem{Ros:NcSch}
{\sc A.~L. Rosenberg}, Noncommutative schemes. \emph{Compositio Math.} \textbf{112}(1998), No. 1, 93--125.

\smallskip \bibitem{Ros:NcSS}
{\sc A.~L. Rosenberg}, {Noncommutative spaces and schemes}.
\emph{Max Planck Institute preprint,}
MPI-1999-84, \emph{Bonn} 1999.

\smallskip \bibitem{Schneider}
{\sc H-J. Schneider}, Principal homogeneous spaces for arbitrary Hopf algebras. Hopf algebras. \emph{Israel J. Math.} \textbf{72}(1990), No. 1-2, 167--195.

\smallskip \bibitem{konradUrsYandl}
{\sc U. Schreiber, C. Schweigert,} and {\sc K. Waldorf},
Unoriented WZW models and holonomy of bundle gerbes. \emph{Comm.
Math. Phys.} \textbf{274}(2007), No. 1, 31--64  \\{\tt
arXiv:hep-th/0512283}.

\smallskip \bibitem{Skoda:coh-states} {\sc  Z.~\v{S}koda},
Coherent states for Hopf algebras. \emph{Lett. Math. Phys.} 
\textbf{81}(2007), No. 1, 1--17  {\tt
arXiv:math.QA/0303357}.

\smallskip \bibitem{Skoda:distr}
{\sc  Z.~\v{S}koda}, {Distributive laws for actions of
monoidal categories}. {\tt arXiv:math.CT/0406310}.

\smallskip \bibitem{Skoda:bicatent}
{\sc  Z.~\v{S}koda}, {Bicategory of entwinings}.
{\tt arXiv:0805.4611}

\smallskip \bibitem{Skoda:lcomp}
{\sc Z.~\v{S}koda}, {Compatibility of (co)actions and localizations}.
{\tt arXiv:0902.1398}.

\smallskip \bibitem{Skoda:Ore}
{\sc Z.~\v{S}koda}, Every quantum minor generates an Ore set.
\emph{Int. Math. Res. Not. IMRN} \textbf{2008}, No. 16, Art. ID
rnn063, 1--8.  {\tt arXiv:math.QA/0604610}.

\smallskip \bibitem{Skoda:eqdis}
{\sc Z. \v{S}koda}, {Equivariant monads and equivariant lifts
versus a 2-category of distributive laws}. {\tt arXiv:0707.1609}.


\smallskip \bibitem{Skoda:ban}
{\sc Z.~\v{S}koda}, Localizations for construction of quantum coset spaces. \emph{Noncommutative geometry and quantum groups} (\emph{Warsaw}, 2001), 265--298, \emph{Banach Center Publ.,} 61, \emph{Polish Acad. Sci., Warsaw,} 2003, {\tt arXiv:math.QA/0301090}.

\smallskip \bibitem{skoda:ncloclec}
{\sc Z.~\v{S}koda}, Noncommutative localization in noncommutative geometry. \emph{Non-commutative localization in algebra and topology,} 220--313, \emph{London Math. Soc. Lecture Note Ser.,} 330, \emph{Cambridge Univ. Press, Cambridge,} 2006,  {\tt math.QA/0403276}.

\smallskip \bibitem{Skoda:qheap}
{\sc Z.~\v{S}koda}, Quantum heaps, cops and heapy categories. \emph{Math. Commun.} \textbf{12}(2007), no. 1, 1--9, {\tt arXiv:math.QA/0701749}.

\smallskip \bibitem{Street:formal}
{\sc R.~Street},
The formal theory of monads. \emph{J. Pure Appl. Algebra} \textbf{2}(1972), No. 2, 149--168;
part II (with {\sc
S.~Lack}) \emph{J. Pure Appl. Algebra} \textbf{175}(2002), No. 1-3, 243--265.


\smallskip \bibitem{Vistoli}
{\sc A. Vistoli},
Grothendieck topologies, fibered categories and descent theory. \emph{Fundamental algebraic geometry,} 1--104, \emph{Math. Surveys Monogr.,} 123, \emph{Amer. Math. Soc., Providence, RI,} 2005, {\tt arXiv:math.AG/0412512}.
\end{thebibliography}
\end{document}